\documentclass{svproc}
\usepackage{amsmath,amsfonts}
\usepackage{mathrsfs}
\usepackage{color,xcolor}
\usepackage{ifpdf}
\usepackage{psfrag}
\usepackage{graphicx,graphics,subfigure,epsfig}
\usepackage{url}
\usepackage{hyperref}
\usepackage{xspace}
\usepackage{nicefrac}

\newcommand{\dpar}[2]{\dfrac{\partial #1}{\partial #2}}

 \newcommand{\R}{\mathbb R}
 \newcommand{\Z}{\mathbb Z}

\renewcommand{\P}{\mathbb P}

\newcommand{\PP}{\mathbb P}

\newcommand{\bbf}{{\mathbf {f}}}
\newcommand{\bn}{{\mathbf {n}}}

\newcommand{\bu}{\mathbf{u}}
\newcommand{\bv}{\mathbf{v}}

\newcommand{\bbu}{\mathbf{u}}
\newcommand{\bbg}{\mathbf{g}}
\newcommand{\hbbg}{\hat{\mathbf{g}}}
\newcommand{\hbbf}{\hat{\mathbf{f}}}

\newcommand{\bx}{\mathbf{x}}

\newcommand{\br}{\mathbf{r}}

\newcommand{\remi}[1]{{#1}}
\usepackage{url}

\begin{document}
\mainmatter              
\title{A personal discussion on conservation, and how to formulate it}
\titlerunning{About conservation}  
%
\author{R\'emi Abgrall \inst{1}}
\authorrunning{R. Abgrall}
%
\tocauthor{R\'emi Abgrall}
\institute{Institute f\"ur Mathematik, Z\"urich Universit\"at, CH 8057 Z\"urich, Switzerland
\email{remi.abgrall@math.uzh.ch},\\ web:
\texttt{\url{https://www.math.uzh.ch/index.php?id=people&key1=8882}}
}
\maketitle              

\begin{abstract}
Since the celebrated theorem of Lax and Wendroff, we know a necessary condition that any  numerical scheme for hyperbolic problem should satisfy: it should be written in flux form. A variant can  also be formulated for the entropy. Even though some schemes, as for example those using continuous finite element, do not formally cast into this framework, it is a very convenient one. In this paper, we revisit this, introduce a different notion of  local conservation which contains the previous one in one space dimension, and explore its consequences. This gives a more flexible framework that allows to get, systematically, entropy stable schemes, entropy dissipative ones, or accomodate more constraints. In particular, we can show that continuous finite element method can be rewritten in the finite volume framework, and all the quantities involved are explicitly computable. We end by presenting the only counter example we are aware of, i.e a scheme that seems not to be rewritten as a finite volume scheme.
\end{abstract}
In this paper, we will consider approximations of the hyperbolic problem
\begin{subequations}
\label{eq:problem}
\begin{equation}
\label{eq:1}
\dpar{\bu}{t}+\text{ div } \bbf(\bu)=0 
\end{equation}
which, for simplicity, we will assume to be defined on $\R^d$ ($d=1,2,3$), and $\bu$ is defined on $\R^d\times [0,T[$, $T>0$, with values in $\Omega\subset \R^p$,  $\Omega$ open. The flux $\bbf=(\bbf_1, \ldots , \bbf_d)$ is such that each $\bbf_j$ are defined on $\Omega$ with values in $\R^p$. The flux $\bbf$ is assumed to be $C^1$ on $\Omega^d$, and the hyperbolicity means that for any $\bn=(n_1, \ldots , n_d)\in \R^d$,  the matrix
$$\nabla \bbf\cdot \bn:=\sum_{j=1}^d \dpar{\bbf_j}{\bbu}n_j$$
is diagonalisable in $\R$.
The PDE \eqref{eq:1} is supplemented with initial condition
$$\bbu(\bx,0)=\bbu_0(\bx)$$
for some initial condition with values in $\Omega$. We will no do any theoretical consideration about this problem, and refer to \cite{dafermos} for more details.
\end{subequations}

In order to integrate \eqref{eq:problem}, several choices need to be done: the type of mesh, the type of functional approximation, how to discretise the divergence term, and finally the time stepping strategy. Often,  one considers a tessaletion of the computational domain by polygons. A common choice is to choose an approximation space $V_h$ which is a subset of $L^2(\R^d)$, then a variational approximation allows to approximate the divergence term, and finally the method of lines is used for time discretisation. Often, the choice $V_h\subset L^2(\R^d)$ hides a common belief that global continuity is not a good idea for \eqref{eq:problem} which allows discontinuous solutions.

In addition, the system \eqref{eq:problem} is complemented by a differential inegality about a state function, the entropy $\eta$, a convex function of $\bu\in \Omega$ (and hence $\Omega$ is assumed to be convex from now on). There exists also $\bbg=(\bbg_1, \ldots , \bbg_d)$, $C^1$ such that for any $j=1, \ldots , d$
$$\nabla_\bu \eta\dpar{\bbf_i}{\bu}=\dpar{\bbg_i}{\bu}.$$
From this, we see that
$$\dpar{\eta}{t}+\text{ div }\bbg =0.$$
For non smooth solution, we impose
\begin{equation}
\label{eq:entropy}
\dpar{\eta}{t}+\text{ div }\bbg \leq 0
\end{equation}
which is true in the sense of distribution.

The existence of \eqref{eq:entropy} is motivated by the canonical example: the Euler equations. Here, the conserved variables are $\bu=(\rho, \rho \bv, E)^T$ where $\rho$ is the density, $\bv$ is the fluid velocity, and $E=e+\frac{1}{2}\rho \bv^2$ is the total energy which is the sum of the internal energy $e$ and the kinetic energy $\frac{1}{2}\rho \bv^2$.
In the simplest case, one can write the internal energy as a function of the density and the pressure $p$, $e=e(\rho, p)$, and the entropy is also a function of these variables, $\eta=\eta(\rho, p)$. We can navigate from one thermodynamic set of two variables, for example $\{\rho, p\}$ to $\{e,\eta\}$ and vice versa. We do enter into more thermodynamic consideration, see eg. \cite{GodelwskiRaviart}.
The conserved variable satisfy \eqref{eq:problem} with the flux defined by
\begin{equation}
\label{eq:euler}
\bbf=\begin{pmatrix}
\rho \bv\\
\rho \bv\otimes\bv+p\text{Id}_{d\times d}\\
(E+p)\bv
\end{pmatrix}
\end{equation}
The entropy flux is
$\bbg= \bv \eta$ 
and in the simplest thermodynamics case, namely the case of a caloricaly perfect gas,
$p=(\gamma-1)e$ where $\gamma$ is a constant and $\eta=\rho \big ( \log p-\gamma\log \rho)-\eta_0$ where $\eta_0$ is a reference.
The system \eqref{eq:problem} with the flux \eqref{eq:euler} is hyperbolic in the domain 
$$\Omega=\{\bu \text{ such that } \rho>0 \text{ and } e>0\},$$ and the entropy is a convex function of $\bu$ if $\gamma>1$.

The simplest example of non linear problem is the Burgers equation (which can be obtained from the Euler equation in the case $\gamma=3$). It is written (in the inviscid case) as
\begin{equation}\label{eq:burger}\dpar{u}{t}+u\dpar{u}{x}=0 \text{ or } \dpar{u}{t}+\frac{1}{2}\dpar{u^2}{x}=0.\end{equation}
Taking $u(x,0)=\sin(\pi x)$, $x\in [-1,1]$, it is easy to see with the method of characteristics that after $t^\star=\tfrac{1}{\pi}$ the solution cannot be smooth, so that the relation \eqref{eq:burger} on the left has no meaning. The same would held for the non conservative form of the Euler equation, for example
\begin{equation}
\label{euler:NC}
\dpar{}{t}\begin{pmatrix}
\rho\\ \bu\\ b\end{pmatrix} + \begin{pmatrix} \text{ div }\big (\rho \bu\big )\\
\big (\bu\cdot \nabla\big ) \bu+\frac{\nabla p}{\rho}\\
\bu\cdot \nabla p+\rho c^2 \text{div }\bu\end{pmatrix}=0
\end{equation}
though these relations are more interesting for practitioners  (since we have a direct access to the pressure and the velocity).

Motivated by this, P. Lax has formalized \remi{what was} already known in the engineering community by taking volumes and looking at what is getting into and out of them.
It is well known that a smooth function is a solution  of \eqref{eq:1} if and only if for any smooth test function $\varphi \in C_0^1(\R^d\times \R^+)$ with compact support, we have
\begin{equation}
\label{eq:weak}
\int_{\R^d\times \R^+} \big (\dpar{\varphi}{t}\bu+\nabla\varphi \cdot \bbf (\bu)\big )\; d\bx\; dt +\int_{\R^d}\varphi(\bx,0)u_0(\bx) \; d\bx=0
\end{equation}
From this we can define the notion of weak solution. Similarly, we have the weak form of the entropy inequality: taking any positive test function, we get
\begin{equation}
\label{eq:weak:entr}
\int_{\R^d\times \R^+} \big (\dpar{\varphi}{t}\eta(\bu)+\nabla\varphi \cdot \bbg (\bu)\big )\; d\bx\; dt +\int_{\R^d}\varphi(\bx,0)\eta_0(\bx) \; d\bx\geq0
\end{equation}
This notion of weak solution is the guiding line of numerical discretisation. The celebrated Lax-Wendroff theorem  is "simply" a way to mimic  this notion, and it provides a generic from of numerical schemes that allows to guaranty, under natural conditions, the convergence to weak entropy solutions. For example, in one dimension, 
\begin{itemize}
\item  being given a regular mesh $\{x_j\}_{j\in \Z}$, and defining control volumes $C_j=(x_{j-1/2}, x_{j+1/2})$ with $x_{j+1/2}=\tfrac{x_j+x_{j+1}}{2}$,
\item  given an initialization\footnote{\remi{If the set $S$ is discrete, $\vert S\vert$ is its cardinal. If $S$ is part of a domain, it is its measure with respect to the lebesgue measure, i.e. its lenght/area/volume.}}
$$u_j^0\approx \frac{1}{\vert C_{j}\vert}\int_{C_{j}} u_0(\bx)\; d\bx$$
\item A numerical scheme of the form
$$\bu_j^{n+1}=\bu_j^n-\frac{\Delta t}{\vert C_{j}\vert}\big ( \hbbf_{j+1/2}-\hbbf_{j-1/2}\big )$$
where $\hbbf_{k+1/2}= \hbbf(u_{k-p}, \ldots, \bu_k\, \ldots , \bu_{k+p})$ is a Lipschitz continuous function depending  $2p+1$ arguments centered around $\bu_k$.
\end{itemize}
Then we know that if the numerical flux is consistant : $\hbbf(u, \ldots , u)=\bbf(u)$, if the the sequence $\{\bu_j^n\}$ is bounded under $\Delta t/\vert C_j\vert \leq C$ and such that one subsequence converges to some $\bv$ in $L^2$, then $\bv$ is a weak solution. The same applies for entropy solutions.

\bigskip
Since about 60 years, all the research has turn around this result and variants. Is this the end of the story? Certainly not.
Several natural questions arise:
\begin{itemize}
\item What happens when a scheme has no longer a flux form ?
Concerning this question, a partial answer was given by Hou and Le Floch in \cite{hou} where they show that a scheme written in incremental form
$$\bu_i^{n+1}=\bu_i^n-C_{i-1/2}^n(\bu_i^n-\bu_{i-1}^n)+D_{i+1/2}^n(\remi{\bu_{i+1}^- - \bu_i^n}),$$ under suitable positivity constraints on the coefficients $C_{l+1/2}$ and $D_{l+1/2}$, a CFL type condition, that a subsequence converges to a function that is a weak solution of
$$\dpar{u}{t}+\dpar{f}{x}=\mu$$ where $\mu$ is  Borel measure. The measure $\mu$ is conjectured to sit on the discontinuities of the solution. Not much more is said, in particular, $\mu$ could possibly vanish.
\item It is know in any text book that a non linear change of variable is in general not permitted. The canonical example is again the Burgers equation. When we set $v=u^3$, the irregular solutions will satisfy
$$\dpar{v}{t}+\frac{3}{4}\dpar{v^3}{x}=0.$$
but the discontinuities will not travel at the same speeds, so that the shocks are different.
However, it would be very interesting to have the possibility to change variable, think, for example,  of the Euler equations in conservative and primitive variables.
\end{itemize}

The purpose of this paper is to explain some ways to overcome this two obstacles. The format of this paper is as follows. We first recall the classical setting of finite volume schemes and the notion of numerical flux, as well as the classical Lax-Wendroff theorem. Then we show a rewriting  of the local conservation condition for these schemes, and describe several classical numerical schemes that satisfy this condition. We show that this condition is then equivalent to the existence of numerical flux, the only difference is that in general these flux are not standard.

Then using this condition, we show several extensions using a non conservative form of a conservative problem, how to modify a scheme to satisfy one or more additional conservative constraints (such an entropy inequality), the use of stagerred grids. We conclude by show an example of scheme that does not seem to be cast in the same framework, though can be shown leading to proper weak solutions of the problem.
\section{Classical conservation versus RD}
\begin{figure}[h]
\begin{center}
\subfigure[]{\includegraphics[width=0.35\textwidth]{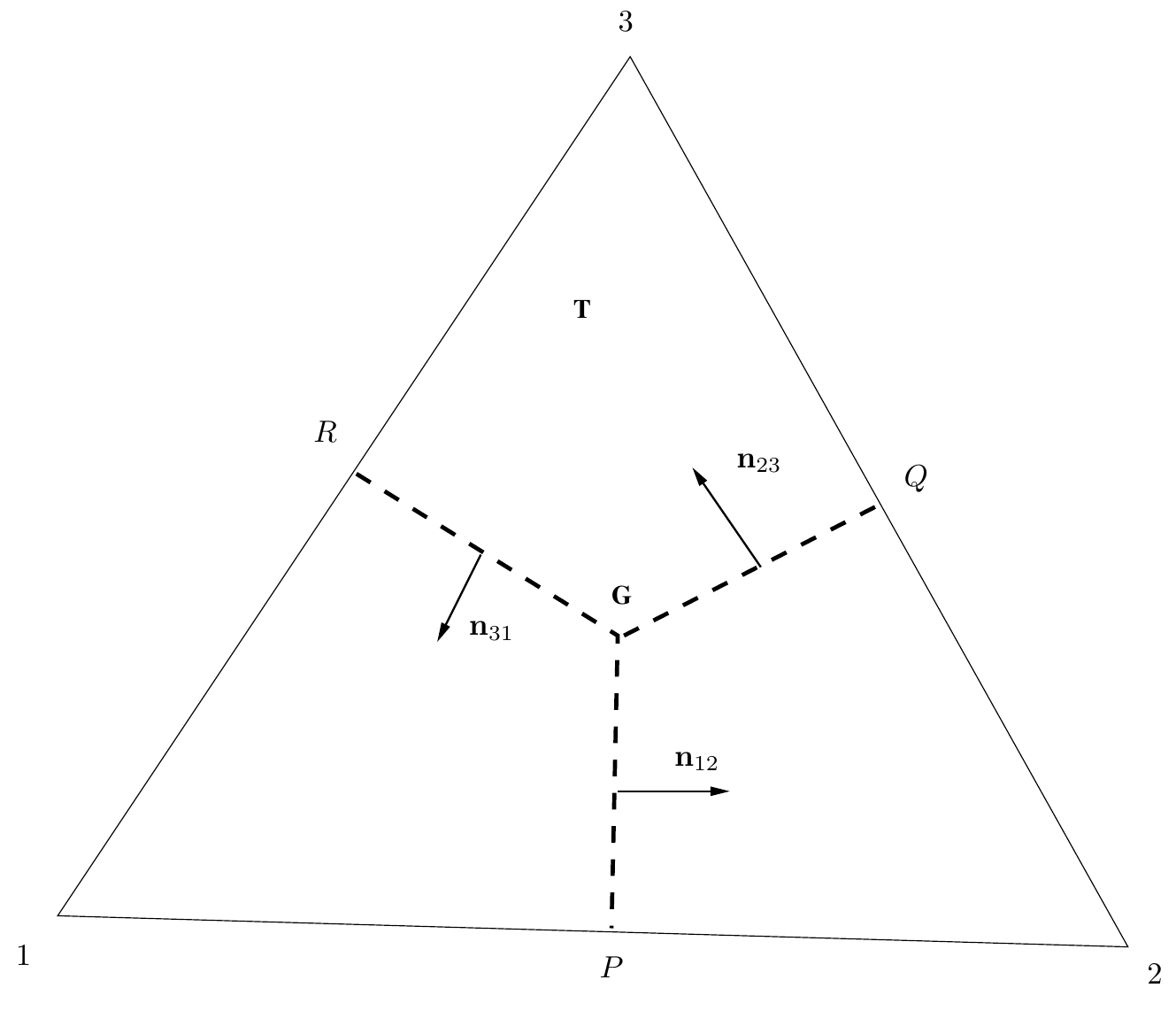}}
\subfigure[]{\includegraphics[width=0.35\textwidth]{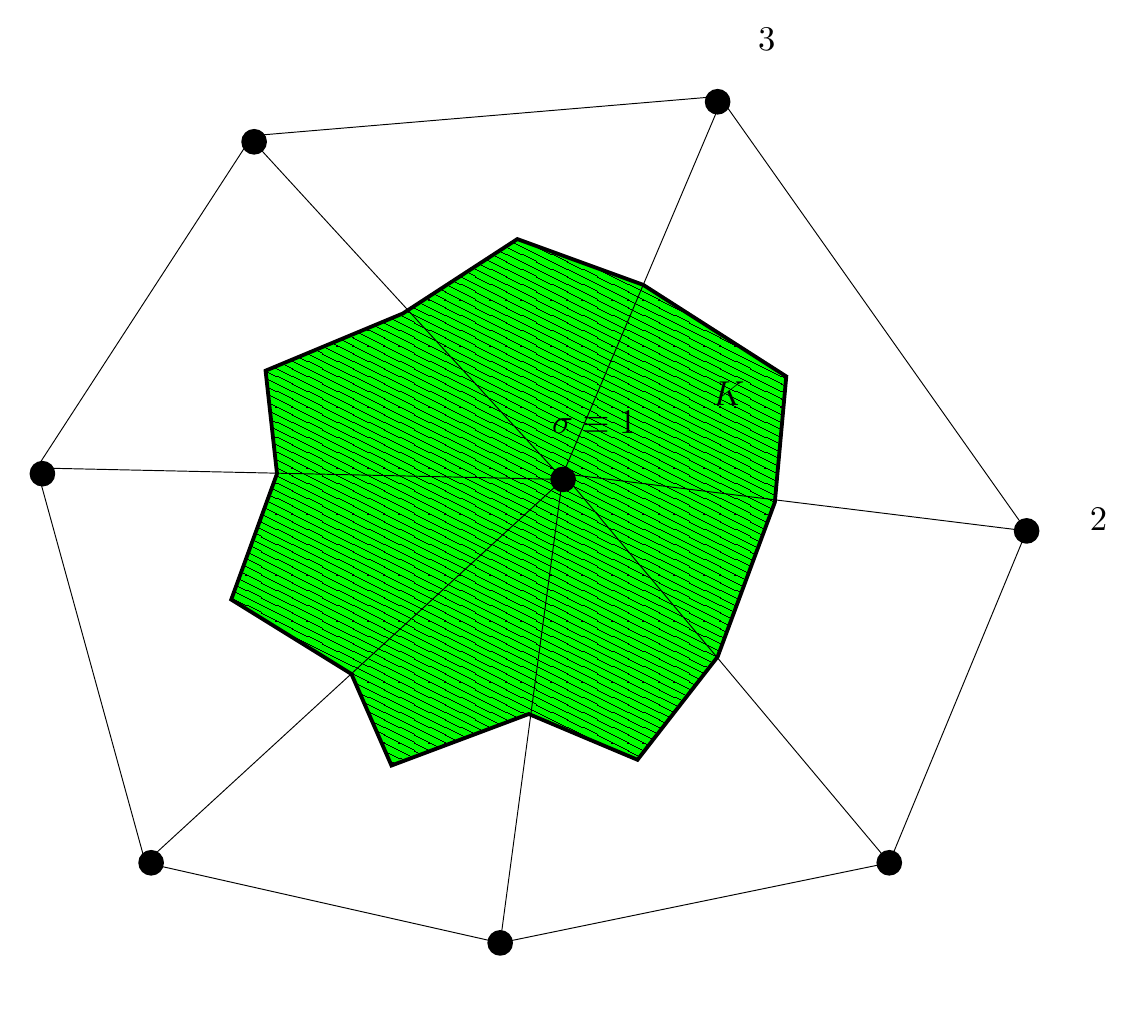}}
\end{center}
\caption{\label{fig:fv} Notations for the finite volume schemes. On the left: definition of the control volume for the degree of freedom $\sigma$.
 The vertex $\sigma$ plays the role of the vertex $1$ on the left picture for the triangle K. The control volume $C_\sigma$ associated to $\sigma=1$ is green on the right and corresponds to $1PGR$ on the left. The vectors $\bn_{ij}$ are normal to the internal edges scaled by the corresponding edge length}
\end{figure}

In this section, we rephrase in part the content of \cite{RADroniou}. We first start from a standard finite volume scheme, and rewrite it in  an equivalent form. In the one dimensional case, we simply say that
$\hbbf_{j+1/2}-\hbbf_{j-1/2}=\hbbf_{j+1/2}-\bbf(u_j)+\bbf(u_j)-\hbbf_{j-1/2}=\Phi_i^{[x_i, x_{i+1}]}+\Phi_i^{[x_{i-1},x_i]}$ with
$$\Phi_i^{[x_i, x_{i+1}]}=\hbbf_{j+1/2}-\bbf(u_j), \quad \Phi_i^{[x_{i-1},x_i]}=\bbf(u_j)-\hbbf_{j-1/2}$$
so that the standard finite volume can be equivalently rewritten as 
\begin{equation}\label{eq:RDS}\bu_i^{n+1}=\bu_i^n-\frac{\Delta t}{\vert C_{i}\vert} \big (\Phi_i^{[x_i, x_{i+1}]}+\Phi_i^{[x_{i-1},x_i]}\big ).\end{equation}
In addition, we note\footnote{This is the essence of Roe's 1981 paper: setting $\Phi_i^{[x_i, x_{i+1}]}=a^-_{i+1/2}(\bu_{i+1}-\bu_i)$ and $\Phi_{i+1}^{[x_i, x_{i+1}]}=a_{j+1/2}^+(\bu_{i+1}-\bu_i)$, we see that the method of characteristics
\eqref{eq:RDS} is conservative if and only if $\bbf(\bu_{i+1})-\bbf(\bu_i)=a_{i+1/2}(\bu_{i+1}-\bu_i)$.}
$$\Phi_i^{[x_i, x_{i+1}]}+\Phi_{i+1}^{[x_i, x_{i+1}]}=\bbf(\bu_{i+1})-\bbf(\bu_i).$$
This can be extended to any kind of finite volume scheme. Instead of going into the full generality, let us take an example to show the principle. For this $\R^d$ is covered by non overlapping simplex denoted by $K$. The vertices of the mesh  are denoted by $\{\sigma_j\}$. For any $\sigma$, we consider a control volume obtained by joining the centroids $\bx_K$ of the simplex sharing $\sigma$ and the mid points of the edges coming out of $\sigma$, see figure \ref{fig:fv} for the notations.
Again, we specialize ourselves to the case of triangular elements, but  \emph{exactly the same arguments} can be given for more general elements,
 provided a conformal approximation space can be constructed. This is  the case for
triangle elements, and we can take $k=1$.

  Since the boundary of $C_\sigma$ is a closed polygon, the scaled outward normals $\bn_\gamma$ to $\partial C_\sigma$ sum up to 0:
$$
\sum_{\gamma \subset \partial C_\sigma}\bn_\gamma=0$$
where $\gamma$ is any of the segment included in $\partial C_\sigma$, such as $PG$ on Figure \ref{fig:fv}.  The finite volume scheme writes
\begin{equation}
\label{eq:FV}\vert C_\sigma\vert \bu_\sigma^{n+1}=\vert C_\sigma\vert \bu_\sigma^{n}-\Delta t\sum_{\gamma \subset \partial C_\sigma} \hbbf_{\bn_\gamma }(\bu_\sigma , \bu^-_\gamma )
\end{equation}
where $\hbbf$ is a consistant numerical flux, and $\bu^-_\gamma$ is the argument on the other side of $\gamma$. It can be evaluated via the MUSCL method, for example. Looking at the spatial increment, we see that
\begin{equation*}
\begin{split}
\sum_{\gamma \subset \partial C_\sigma} \hbbf_{\bn_\gamma }(\bu_\sigma , \bu^-& )= \sum_{\gamma \subset \partial C_\sigma} \hbbf_{\bn_\gamma }(\bu_\sigma, \bu^- )- \bigg (\sum_{\gamma \subset \partial C_\sigma}\bn_\gamma\bigg )\cdot \bbf (\bu_\sigma)\\
&=\sum\limits_{K, \sigma\in K} \sum\limits_{\gamma \subset \partial C_\sigma\cap K} \big ( \hbbf_{\bn_\gamma }(\bu_\sigma, \bu^- )-\bbf (\bu_\sigma)\cdot \bn_\gamma \big )
\end{split}
\end{equation*}
To make things explicit, in $K$, the internal boundaries are $PG$, $QG$ and $RG$, and those around $\sigma\equiv 1$ are $PG$ and $RG$.
We set
\begin{equation}
\begin{split}
\Phi_\sigma^K(\bu^h)&=\sum\limits_{\gamma\subset \partial C_\sigma\cap K} \big ( \hbbf_{\bn_\gamma }(\bu_\sigma, \bu^- )-\bbf (\bu_\sigma)\cdot \bn_\gamma \big )\\
&=\sum\limits_{\gamma\subset \partial ( C_\sigma\cap K )}  \hbbf_{\bn_\gamma }(\bu_\sigma, \bu^- ).
\end{split}
\label{fv:res:sigma}
\end{equation}
The last relation uses the consistency of the flux and the fact that $C_\sigma\cap K$ is a closed polygon. The quantity $\Phi_\sigma^K(\bu^h)$ is the normal flux on $C_\sigma\cap K$.
\remi{If now we sum up these three quantities, we get:}
\begin{equation*}
\begin{split}
\sum_{\sigma\in K} \Phi_\sigma^K(\bu_h)&= \bigg ( \hbbf_{\bn_{12}}(\bu_1,\bu_2)-\hbbf_{\bn_{13}}(\bu_1,\bu_3)-\bbf(\bu_1)\cdot\bn_{12}+\bbf(\bu_1)\cdot\bn_{31}\bigg )\\
&+\bigg ( \hbbf_{\bn_{23}}(\bu_2,\bu_3)-\hbbf_{\bn_{12}}(\bu_2,\bu_1)+\bbf(\bu_2)\cdot\bn_{12}-\bbf(\bu_2)\cdot\bn_{23}\bigg )\\
&+\bigg ( -\hbbf_{\bn_{23}}(\bu_3,\bu_2)+\hbbf_{\bn_{31}}(\bu_3,\bu_1)-\bbf(\bu_3)\cdot\bn_{23}+\bbf(\bu_3)\cdot\bn_{31}\bigg )\\
&= \bbf(\bu_1)\cdot \big ( \bn_{12}-\bn_{31}\big ) +\bbf(\bu_2)\cdot \big ( -\bn_{23}+\bn_{31}\big )
+\bbf(\bu_3)\cdot \big ( \bn_{31}-\bn_{23}\big )\\
&=\bbf(\bu_1)\cdot\frac{\bn_1}{2}+\bbf(\bu_2)\cdot\frac{\bn_2}{2}+\bbf(\bu_3)\cdot\frac{\bn_3}{2}
\end{split}
\end{equation*}
where $\bn_j$ is the scaled inward normal of the edge opposite to vertex $\sigma_j$, i.e. twice the gradient of the $\PP^1$ basis function
 $\varphi_{\sigma_j}$ associated to this degree of freedom.
Thus, we can reinterpret the sum as the boundary integral of the Lagrange interpolant of the flux.
The finite volume scheme is then a residual distribution scheme with residual defined by \eqref{fv:res:sigma}
and a total residual defined by
\begin{equation}
\label{fv:tot:residu}
\Phi^K:=\int_{\partial K} \bbf^h\cdot \bn , \qquad \bbf^h=\sum_{\sigma\in K} \bbf(\bu_\sigma)\varphi_\sigma.
\end{equation}

Form now on, we will assume that the domain can be split into polygons $K$ (above it was simplex) which vertex  will be the $\sigma$s. From  these polygones, we  can construct of control volumes noted $C_\sigma$. The other situation is the converse: from a family of control volumes, we can construct polygons $K$ with vertices $\sigma$: the difference is between vertex centered or volume centered schemes in the Finite volume vocabulary. We will focus on schemes that can be written in the following form:
\begin{subequations}\label{eq:RD}
\begin{equation}
\label{eq:RD:update}
\vert C_\sigma\vert \bu_\sigma^{n+1}=\vert C_\sigma\vert \bu_\sigma^{n}-\Delta t \sum_{K, \sigma\in K}\Phi_\sigma^K(\bu^n)
\end{equation}
with 
\begin{equation}
\label{eq:RD:conservation}\sum_{\sigma\in K}\Phi_\sigma^K(\bu^n)=\int_{\partial K}\hbbf_\bn\ d\gamma.
\end{equation}
\end{subequations}
In \eqref{eq:RD:update}, $\Phi_\sigma^K(\bu^n)$ is a function that depends on a finite number values of the $\bu_{\sigma'}^n$, it  is assumed to be Lipschitz continuous. In \eqref{eq:RD:conservation} we assume that for any edge/face that is the intersection of two polygons, $\gamma=K\cap K'$,
$$\int_{\gamma\subset K}\hbbf_\bn\ d\gamma+\int_{\gamma\subset K'}\hbbf_\bn\ d\gamma=0.$$
In other words, the flux are the same, up to the sign.

\bigskip

This framework is not adapted uniquely  to finite volume. We can consider, using continuous finite elements, the SUPG scheme. Considering a mesh made of simplex, we take
$$V^h=\{ \bv^h\in C^0(\R^d), \text{ for any }K,  \bv^h_{\vert K} \in \P^r\big\}$$
and look for $\bu^h\in \big(V^h\big )^p$ such that for any $\bv^h\in V^h$ we have
\begin{equation}
\label{SUPG}
\sum\limits_{K\subset \R^d}\bigg ( -\int_K\nabla v^h \cdot \bbf(\bu^h)\; d\bx+\int_{\partial K} v^h\bbf(\bu^h)\cdot \bn \; d\gamma+
h_K\int_K\big ( \nabla_\bu\bbf \cdot \nabla v^h\big ) \tau \big ( \nabla_\bu\bbf\cdot \nabla \bu^h\big ) \;d\bx \bigg )=0
\end{equation}
By continuity, of course the boundary term cancel, but we have written the scheme in this way to exhibit the residuals: if $\{\varphi_\sigma\}$ is the set of Lagrange basis functions,
\begin{equation}
\label{SUPG:residual}
\Phi_\sigma^K= -\int_K\nabla v^h \cdot \bbf(\bu^h)\; d\bx+\int_{\partial K} v^h\bbf(\bu^h)\cdot \bn \; d\gamma+
h_K\int_K\big ( \nabla_\bu\bbf \cdot \nabla v^h\big ) \tau \big ( \nabla_\bu\bbf\cdot \nabla \bu^h\big ) \;d\bx,
\end{equation}
and we get the conservation relation:
\begin{equation}
\label{SUPG:residual:conservation}
\sum_{\sigma\in K}\Phi_\sigma^K=\int_{\partial K} \bbf(\bu^h)\cdot \bn \; d\gamma
\end{equation}
because $\sum\limits_{\sigma\in K}\big (\varphi_\sigma\big )_{|K}=1.$

All the schemes we know, exept maybe one that we will sketch at the end of this paper, can be rewritten in a distribution form: this is true from finite volume (or any order), to dG, via continuous finite element with different stabilisation mechanisms. Some details can be found in \cite{RADroniou}. This can also apply to schemes adapted to the Lagrangian formalism, see \cite{maire} for example.

Let us examine now the converse: assuming a scheme of the form \eqref{eq:RD:update} with \eqref{eq:RD:conservation}, can we identify "flux" so that the scheme \eqref{eq:RD:update} can be rewritten equivalently in the form \eqref{eq:FV}.  First, what is a flux?
\begin{definition}[consistent flux]
A flux function $\hbbf$ is a function that depends on a normal $\bn$ and a set of arguments $\{\bu_1, \ldots, \bu_N\}$ such that:
\begin{enumerate}
\item \remi{$\hbbf$ is continuous with respect to its arguments,}
\item $\hbbf(-\bn; \bu_1, \ldots, \bu_N)=-\hbbf(\bn; \bu_1, \ldots, \bu_N)$
\item It is consistent if for any $\bu$ and $\bn$, $\hbbf(-\bn; \bu, \ldots, \bu)=\bbf(\bu)\cdot \bn$
\end{enumerate}
We will also use the notation $\hbbf_\bn(\bu_1, \ldots, \bu_N)$ or simply $\hbbf_\bn$.
\end{definition}

Any $K$ appearing in the sum \eqref{eq:RD:update}, contains a set of degrees of freedom, say $\mathcal{S}=\{\sigma_1, \ldots , \sigma_m\}$ from which we can construct a graph, or a triangulation, which vertex are the $\sigma_i$s. An important property of this graph is that it is simply connected.
Then the question is to find $\{\hbbf_{\sigma\sigma'}\}_{\sigma, \sigma'\in \mathcal{S}}$ such that
\begin{subequations}\label{GC:1}
\begin{equation}
\label{GC:1.1}
\Phi_\sigma=\sum_{\text{ edges }[\sigma,\sigma']} \hbbf_{\sigma,\sigma'}+\hbbf_\sigma^{b}
\end{equation}
with
\begin{equation}
\hbbf_{\sigma,\sigma'}=-\hbbf_{\sigma',\sigma}
\label{GC:1.2}
\end{equation}
and $\hbbf_\sigma^{b}$ is the 'part' of $\int_{\partial K} \hbbf_\bn(\bu^h,\bu^{h,-}) \; d\gamma$ associated to $\sigma$. {The control volumes will be defined by their normals so that we get consistency.} The normal will be defined later, as well as the control volumes.

Note that \eqref{GC:1.2} implies the conservation relation
\begin{equation}
\label{GC:conservation}
\sum\limits_{\sigma\in K}\Phi_\sigma=\sum\limits_{\sigma\in K}\hbbf_\sigma^b.
\end{equation}
In short, we will take
\begin{equation}
\label{BC:1.3}
\hbbf_\sigma^b=\oint_{\partial K} \varphi_\sigma\; \hbbf_\bn (\bu^h,\bu^{h,-}) \; d\gamma,
\end{equation}
\end{subequations}
but other  examples can be considered provided the consistency \eqref{GC:conservation} relation holds true, see \cite{RADroniou}.
Any edge $[\sigma,\sigma']$ is either direct or, if not, $[\sigma',\sigma]$ is direct. Because of \eqref{GC:1.2}, we only need to know $\hbbf_{\sigma,\sigma'}$ for direct edges. Thus we introduce the notation $\hbbf_{\{\sigma,\sigma'\}}$ for  the flux  assigned to  the direct edge whose extremities are $\sigma$ and $\sigma'$. We can rewrite \eqref{GC:1.1} as, for any $\sigma\in \mathcal{S}$,
\begin{equation}
\label{GC:1.1bis}
\sum_{\sigma'\in \mathcal{S}} \varepsilon_{\sigma,\sigma'} \hbbf_{\{\sigma,\sigma'\}}=\Psi_\sigma:=\Phi_\sigma-\hbbf_\sigma^b,
\end{equation}
with $$
\varepsilon_{\sigma,\sigma'}=\left \{
\begin{array}{ll}
0& \text{ if }\sigma \text{ and }\sigma' \text{ are not on the same edge of }\mathcal{T},\\
1& \text{ if } [\sigma,\sigma']\text{ is an edge and } \sigma \rightarrow \sigma' \text{ is direct,}\\
-1&  \text{ if } [\sigma,\sigma']\text{ is an edge and } \sigma' \rightarrow \sigma \text{ is direct.}
\end{array}
\right .
$$

Hence the problem is to find  a vector $\hbbf=(\hbbf_{\{\sigma,\sigma'\}})_{\{\sigma,\sigma'\} \text{ direct edges}}$ such that
$$A\hbbf=\Psi$$
where $\Psi=(\Psi_\sigma)_{\sigma\in \mathcal{S}}$ and $A_{\sigma \sigma'}=\varepsilon_{\sigma,\sigma'}$.

We have  the following lemma \cite{RADroniou}which shows the existence of a solution.
\begin{lemma}\label{lemma:flux}
For any couple $\{\Phi_\sigma\}_{\sigma\in \mathcal{S}}$ and $\{\hbbf_\sigma^{b}\}_{\sigma\in \mathcal{S}}$ satisfying the condition  \eqref{GC:conservation}, there exists numerical flux functions $\hbbf_{\sigma,\sigma'}$ that satisfy \eqref{GC:1}. Recalling that the  matrix of the Laplacian of the graph is $L=AA^T$, we have
\begin{enumerate}
\item The rank of $L$ is $\mathcal{S}|-1$ and its image is $\big (\text{span}\{\mathbf{1}\})^\bot$. We still denote the inverse of $L$ on $\big (\text{span}\{\mathbf{1}\} )^\bot$ by $L^{-1}$,
\item 
With the previous notations, a solution is 
\begin{equation}
\label{eq:lemma}\big (\hbbf_{\{\sigma,\sigma'\}}\big )_{\{\sigma,\sigma'\} \text{ direct edges}}=A^TL^{-1} \big (\Psi_\sigma\big )_{\sigma\in \mathcal{S}}.\end{equation}
\end{enumerate}
\end{lemma}

\begin{figure}[htb]
\begin{center}
\subfigure[]{\includegraphics[width=0.35\textwidth]{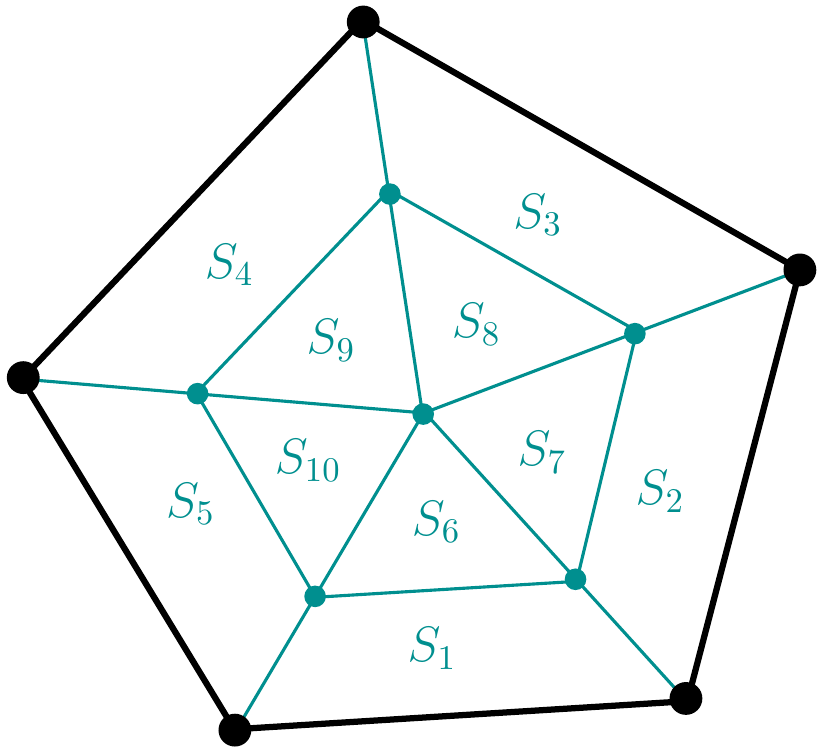}}
\subfigure[]{\includegraphics[width=0.35\textwidth]{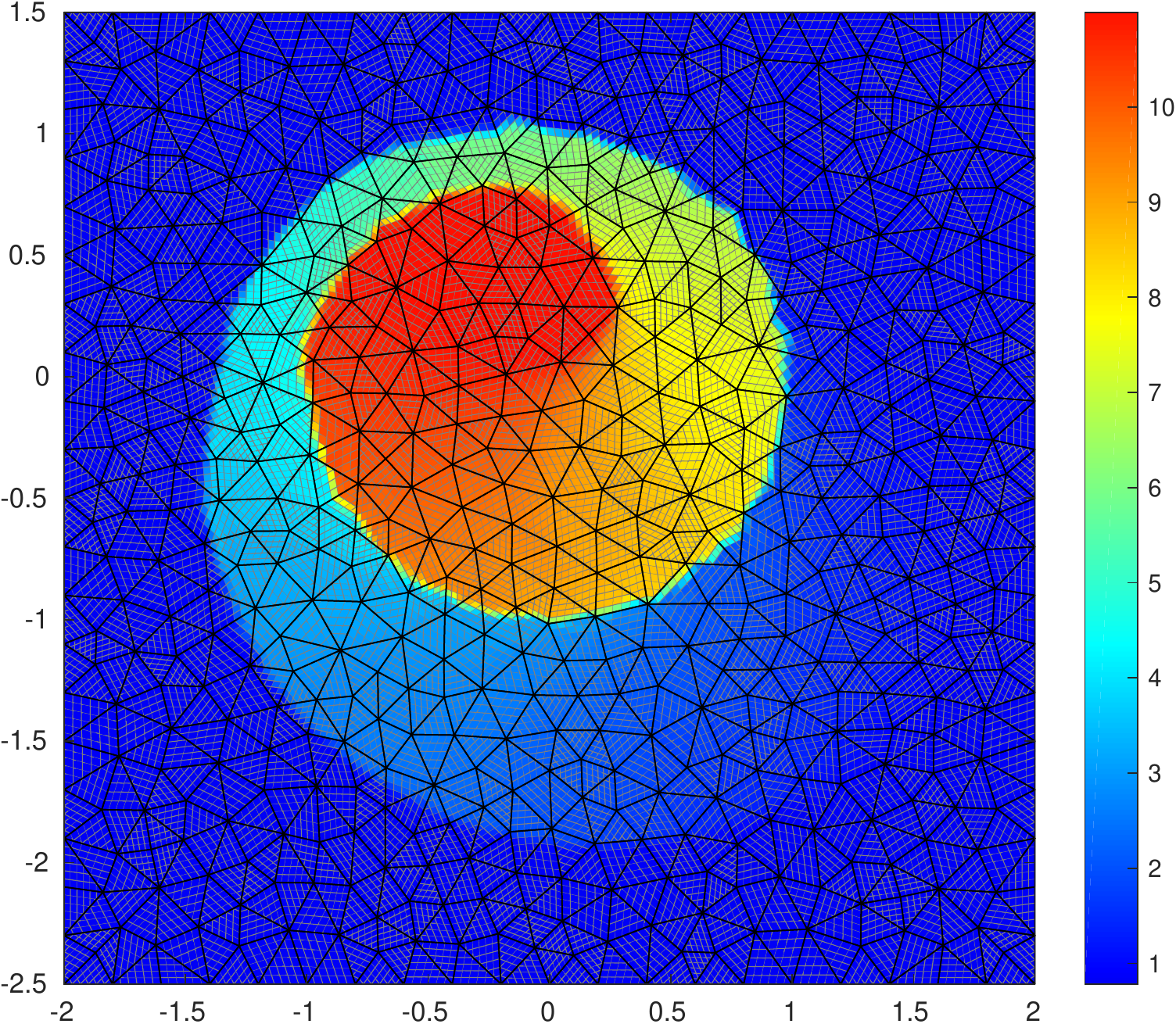}}
\end{center}
\caption{\label{francois} (a): example of a subdivision, (b) solution obtained with the method of \cite{vilar} for the KPP problem. These figures have been generated by Fran\c{c}ois Vilar, Universit\'e de Montpellier.}
\end{figure}

This result has been used to develop in \cite{vilar} a numerical scheme  of dG type, where, when a criteria explaining that some problem occur (negative density, or pressure, or creation of artificial oscillation), the elements are splitted into  sub elements where a low order finite volume scheme  is applied. An example of sub-cells and application if show in figure \ref{francois}.

A convergence proof (in the statistical solution framework) of the schemes \eqref{eq:RD:update} with \eqref{eq:RD:conservation} can be found in \cite{Oeffner}.

The reinterpretation of schemes in term or residuals has other applications. For example in \cite{Paola}, we have considered a mixture of perfect gases. It is well known that because of the possible incompatibility between the numerical dissipations attached to each of the conserved variables, the pressure and the velocity may oscillate across a contact line. There are several solutions to cure that. One is to start from the Euler equation in primitive variables (as \eqref{euler:NC} but we need to consider two masses), but then one has to do something very specific to take into account the shocks properly. For example, in \cite{karni}, a standard finite volume method is used anywhere, but around the slip line which location is estimated by a level set. In \cite{Paola}, we start directly from the non conservative formulation where the variables are $\rho$, $\bv$ and $e$ the internal energy. If $\Delta$ represents the time increment, we notice that
\begin{equation}
\label{energyupdate}\Delta E=\Delta e+\frac{\bv^{n+1}+\bv^n}{2}\Delta \big (\rho \bv\big )-\frac{\bv^{n+1}\bv^n}{2}\Delta \rho.\end{equation}
Then, for a scheme of the form \eqref{eq:RD:update}, where $\bu=(\rho, \rho \bu, e)$, it is shown in \cite{Paola} that if the residual on the energy (and only this one) is modified such that
\begin{equation}
\label{RD:2fluids}\int_{\partial K} \bbf^E\cdot \bn \; d\gamma =\sum_{\sigma\in K} \Phi_\sigma^e+\sum_{\sigma\in K}\frac{\bv_\sigma^{n+1}+\bv_\sigma^n}{2}\Phi_\sigma^{\rho \bv}+\sum_{\sigma\in K} \frac{\bv_\sigma^{n+1}\bv_\sigma^n}{2}\Phi_\sigma^{\rho},
\end{equation} then the scheme will be locally conservative. The modification is done by adding to the initial energy residuals $\Phi_\sigma^E$ the same quantify $r^K$ for all the degrees of freedom in $K$ such that \eqref{RD:2fluids} holds true. This formulation can be further refined so that one keeps the good behavior of the contact lines, see \cite{Paola} for details and applications with very non linear equations of state. 
\section{Staggering}
In another application, one considers an approximation of \eqref{euler:NC} where the velocity is globaly continuous and the thermodynamic parameters (the density, the pressure, the internal energy) are only in $L^2$. This is an Eulerian version of \cite{svetlana1,svetlana2} which are inspired from \cite{dobrev} which is a finite element generalisation of the Wilkins scheme \cite{wilkins}. The Wilkins scheme  is very popular in the Lagrange hydrodynamics community because use the specific internal energy as a variable though producing the correct weak solutions of the problem.

Here we summarise \cite{staggered}. 
We start from a scheme having the from \eqref{eq:RD:update}. The velocity is approximated by a piecewise polynomial function or degree $r\geq 1$ that is globally continuous. For technical reasons, we assume that we use a Bernstein basis (each of the basis function is positive). The pressure, internal energy and the density are approximated by polynomials of degree $r-1$, and again we use Bernstein basis functions.  The thermodynamical degrees of freedoms are denoted by $\sigma_T$ and the velocity degrees of freedom are $\sigma_V$. The density is in $L^2$, so we use a discontinuous Galerkin approximation, with a Riemann solver (that is used only for the density).
Then, for each degree of freedom,  we discretize the velocity equation in the form
$$\bv_\sigma^{n+1}=\bv_\sigma^n-\frac{\Delta t}{\vert C_\sigma\vert }\sum_{K, \sigma\in K} \Phi_\sigma^K.$$
Here, $\vert C_\sigma\vert$ is the mass of $\varphi_\sigma$ which is positive. We take
$$\Phi_\sigma^{\bv, K}=\int_K \varphi_\sigma \bigg (\big (\bv\cdot \nabla \big )\bv +\frac{\nabla p}{\rho}\big ) \; d\bx$$
computed by numerical quadrature, and 
$$\Phi_\sigma^{e,K}=\int_K\varphi_\sigma\bigg ( \bv\cdot\nabla e+(e+p)\text{ div }\bu\bigg )\; d\bx,
$$ again by numerical quadrature. As such, there is no hope to obtain a good scheme. 

In order to modify it, we start from an inspection of what would be needed to get a Lax-Wendroff like theorem. It turns out that for the \remi{update in time} of the momentum, we have the relation:
$$\rho^{n+1}\bv^{n+1}-\rho^n\bv^n=(\rho^{n+1}-\rho^n)\bv^n+\rho^{n+1}(\bv^{n1}-\bv^n),$$
from which we can infer, after some calculations, that if the mass and velocity residual satisfy, in each element,
\begin{equation}
\label{staggered:vitesse}
\int_{\partial K} \rho^{n}\bv^n\cdot \bn=\sum_{\sigma_V\in K} \theta^\bv_{\sigma_V}\Phi^\bv_{\sigma_V}+\sum_{\sigma_T\in K} \mu_{\sigma_T}^\bv \Phi_{\sigma_T}^\rho
\end{equation}
where the coefficients $\theta^\bv_{\sigma_V}$ and $\mu_{\sigma_T}$ are explicitly computable, see \cite{staggered}. In addition $ \theta^\bv_{\sigma_V}>0$. Of course the initial scheme does not satisfy this constraint, but keeping the density at $t_n$ and $t_{n+1}$ unchanged, we can modify $\Phi^\bv_{\sigma_V}$ by $\Phi^\bv_{\sigma_V}+r^K$, so that \eqref{staggered:vitesse} is true.

\begin{figure}[ht]
\begin{center}
\subfigure[]{\includegraphics[width=0.35\textwidth,clip=]{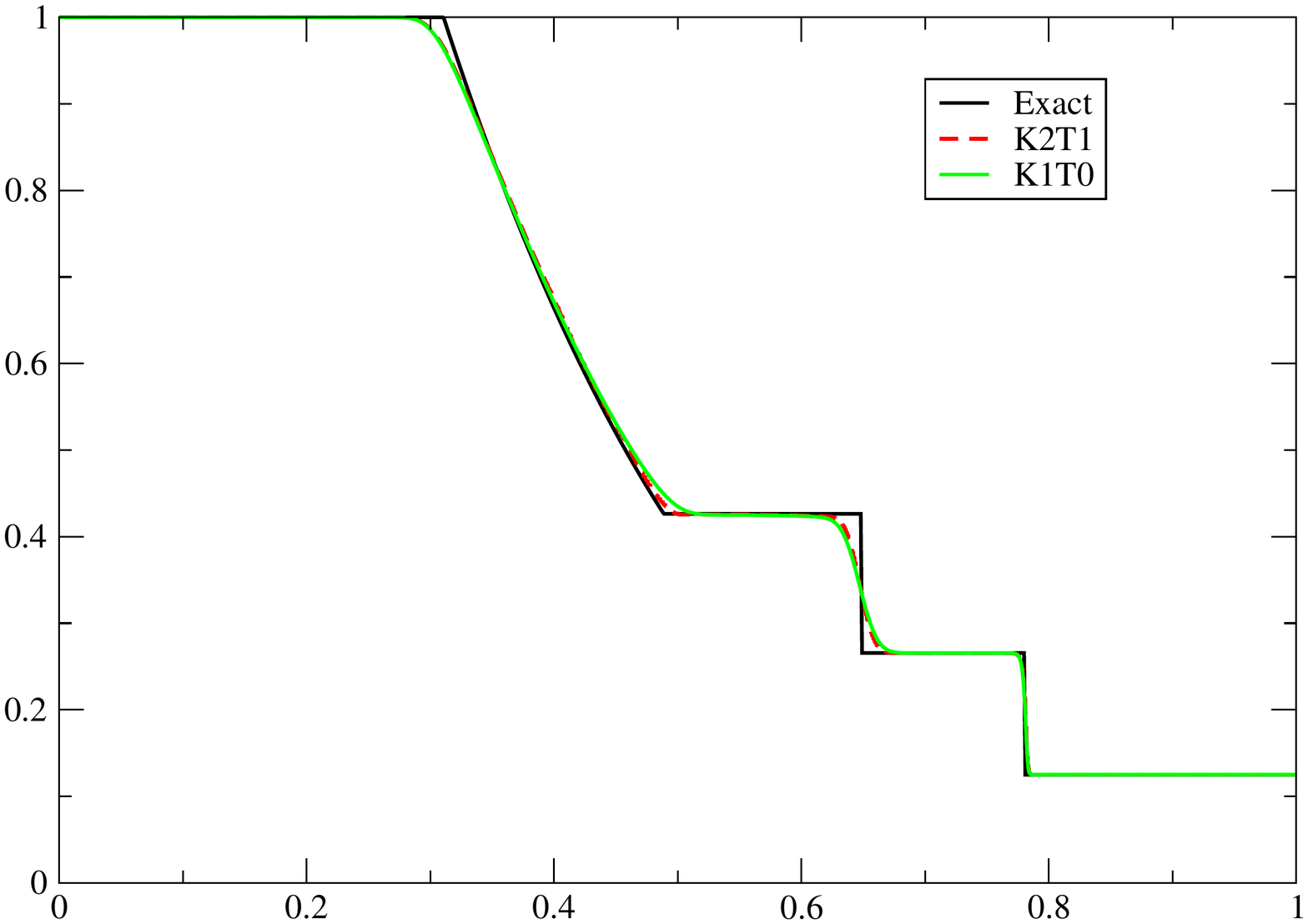}}\subfigure[]{\includegraphics[width=0.35\textwidth,clip=]{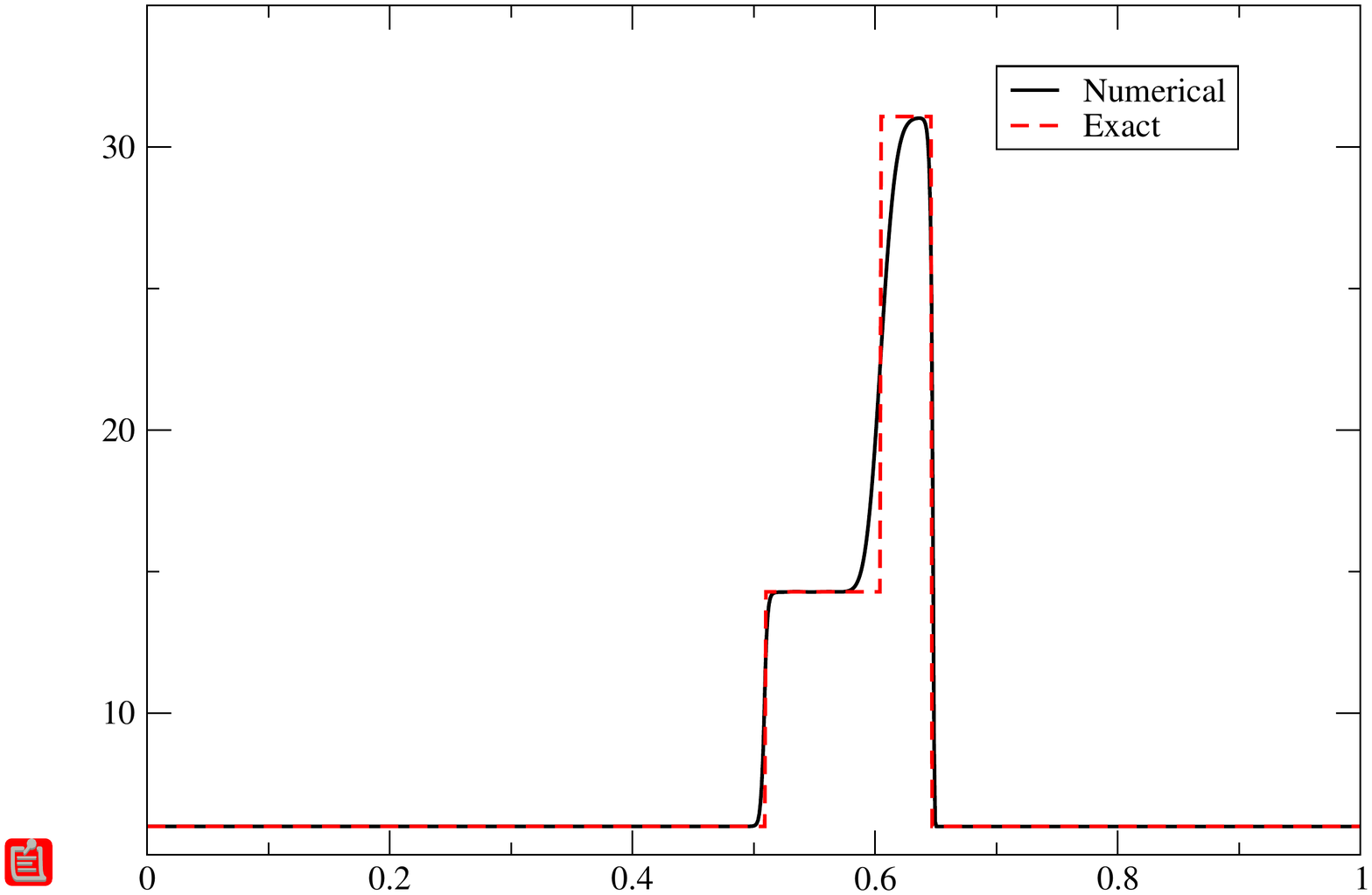}}
\subfigure[]{{\includegraphics[width=0.25\textwidth,clip=]{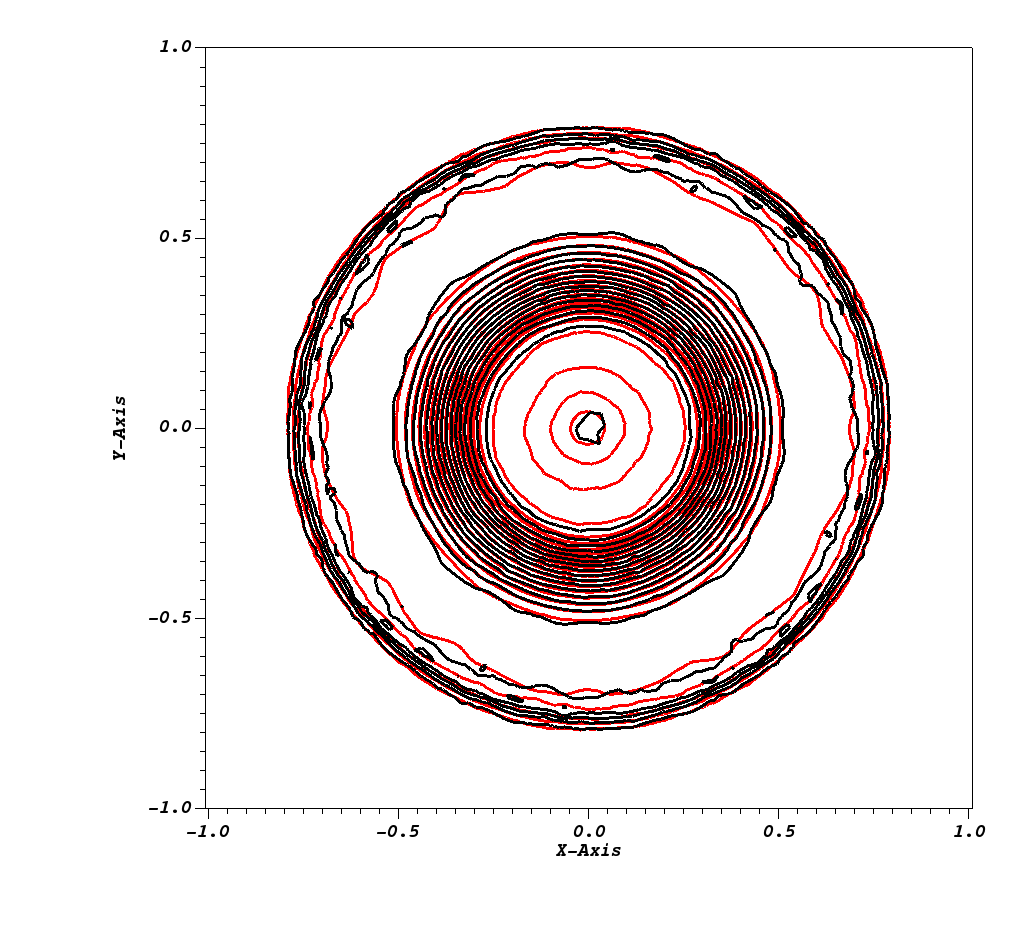}}}
\end{center}
\caption{\label{stagered 1d} (a): density for the Sod test case. We observed the solution without correction (red) ans the numerical one with 1000 points, it superimpose the exact solution. (b) Exact and numerical solutions of the Collela and Woodwards test case for (a) density
at time $T = 0.012$ with $CFL = 0.1$ computed on a mesh with $1000$ points. (c)Comparison of the solutions of the 2D shock test case for the pressure obtained by a conservative scheme (red) and the staggered one (black) at time $T = 0.16$.}
\end{figure}

Similar algebraic manipulation can be done for the energy, and here we need to mimick \eqref{energyupdate}, and this is always possible.
We illustrate by some example taken from \cite{staggered} in figure \ref{stagered 1d}.

\section{Additional conservation laws}
It is also possible to use this framework in order to take into account additional conservation laws, at least in the semi-discrete sense. This part is a short summary of \cite{RAent} and then \cite{RANord1,RANord2,RAPORanocha,KinticMoment}.

Let us consider the entropy, for example. Multiplying \eqref{eq:1} by $\bv:=\nabla_\bu\eta$ we obtain \eqref{eq:entropy}. Hence, considering a scheme of the form
$$
\vert C_\sigma\vert \dfrac{d\bu_\sigma}{dt}+\sum_{K, \sigma\in K}\Phi_\sigma^K=0,$$
we would have
$$\vert C_\sigma\vert \dfrac{d\eta}{dt}+\sum_{K, \sigma\in K} \bv_\sigma\cdot \Phi_\sigma^K=0.$$
This gives the idea of introducing the entropy residuals,
$$\Psi_\sigma^K=\bv_\sigma \cdot \Phi_\sigma^K.$$
However, there is no reason why if 
$$\sum_{\sigma\in K} \Phi_\sigma^K=\int_{\partial K}\hbbf_\bn\; d\gamma$$ we would have $\sum_{\sigma\in K} \Psi_\sigma^K$ related to some boundary integral of the entropy flux.

The trick is similar as before, we introduce $\br_\sigma^K$ and consider $\widetilde{\Phi_\sigma^K}=\Phi_\sigma^K+\br_\sigma^K$ such that we sill have the conservation relation and have in addition
$$\sum_{\sigma\in K} \Psi_\sigma^K\geq \int_{\partial K} \hbbg_\bn\; d\gamma$$ where $\hbbg_\bn$ is some consistant numerical approximation of the entropy flux.
The conservation requirement implies that
$$\sum_{\sigma\in K}\br_\sigma^K=0.$$
This condition can be met if $$\br_\sigma^K=\alpha_K \big ( \bv_\sigma-\overline\bv), \quad \overline\bv=\frac{1}{\#K}\sum_{\sigma\in K}\bv_\sigma,$$
and we get
$$\alpha_K \sum_{\sigma\in K} \big ( \bv_\sigma-\overline\bv)^2\geq \int_{\partial K} \hbbg_\bn\; d\gamma-\sum_{\sigma\in K}\bv_\sigma\cdot \Phi_\sigma^K.$$
In \cite{RAent} \remi{there is a discussion} on the choice of $\hbbf$ such that 
$$\int_{\partial K} \hbbg_\bn\; d\gamma-\sum_{\sigma\in K}\bv_\sigma\cdot \Phi_\sigma^K=O\big ( \bv-\overline\bv)^2$$
to guaranty that $\alpha_K$ does not blow up. The discussion is certainly not finished. 

In \cite{RAPORanocha}, it is shown how to extend this approach when we have several constraints: in that paper, we had discussed the case of the entropy and the kinetic energy preservation. Instead of a simple linear equation, one gets a system of size the number of constraints, and it can be solved by least square. Interestingly, the  more degrees of freedom (i.e. the higher the formal order is), the more constraint one can a priori satisfy. We are not able to prove that if one starts with a stable scheme, the modified scheme will also be stable. \remi{In  practice, we have never observed any instability.}

 \remi{In \cite{RANord1}, using this technique, we have shown that a \remi{fully} centered scheme can be made stable! There is no contradiction with the classical analysis that uses periodicity. Of course this cannot apply to periodic problems, but to problems with proper boundaries. \remi{The idea is, for a transport problem, to write a scheme. In \cite{RANord1}, the examples were considering triangular type  meshes, the unknown was approximated  with $\P^k$ spatial approximation, and the scheme was a simple weak formulation with accurate enough quadrature formula}. In order to have an energy bound, we modify the scheme, so that on each element, we have an exact energy production. This is done following the above ideas. When summing all the contributions,  the interior of the domain does not produce energy, and all has to be controlled at the boundary. In fact, up to this energy summation argument that can be made possible here at the full discrete level, this very classical: it is well know that a dG method is $L^2$ stable if one has a dissipative boundary flux. Hence we are considering a dG method with one element (the whole domain), nobody has ever said that the approximation must be polynomial: it simply need to be accessible to some form of the divergence theorem. All is made in purpose here for that. This also applies to non linear problems, see \cite{RANord2}.}
 
 This technique can be further used. In \cite{KinticMoment}, a similar approach is used so that a local conservation equation for the kinetic momentum is also satisfied. In \cite{RADumbser}, it has been used to obtain a thermodynamicaly compatible scheme for a system that can simultaneously describe a fluid and a solid. There is no space here to describe the method, we refer to the publication. Let us only mention that the variables contains the entropy, not the total energy, and we nevertheless have a method that is locally conservative.
 
\section{A scheme that does not fit in this framework}
In \cite{AF1,AF2}, and following an idea of P.L. Roe \cite{Roe}, we consider a grid with points $x_i<x_{i+1}$. The conserved variable $\bu$ is approximated by its average in each volume $[x_i,x_{i+1}]$, and we also consider the \emph{point} values of some other set of variables denoted by $\bv_i$ (for the point $x_i$, $\bv$ can be $\bu$, but also any transformation of $\bu$ by some regular mapping $\Psi$) and it satisfies the evolution equation
$$\dpar{\bv}{t}+\underbrace{(\nabla_\bu\psi)^{-1} \nabla_\bu \bbf\nabla_\bu \psi}_{J} \; \dpar{\bv}{x}=0.$$
For fluid mechanics, $\bv$ can be the primitive variables, for example.
 In each cell $[x_i,x_{i+1}]$ one has $\overline{\bu}_{i+1/2}$ and the point values $\bv_i=\Psi(\bu_i)$, $\bv_{i+1}=\Psi(\bu_{i+1})$. From this one can construct   a globally continuous reconstruction which is quadratic  in each cell. It amounts to Simpson's formula:
 $$\overline{\bu}_{i+1/2}=\frac{1}{4}\big ( \psi^{-1}(\bv_i)+\psi^{-1}(\bv_{i+1/2})+\psi^{-1}(\bv_{i+1})\big )$$
 so that one can obtain a third order approximation of $\bv_{i+1/2}\approx\bv(\frac{x_i+x_{i+1}}{2}).$
 
 The averaged values and point value  are evolved by 
 $$(x_{i+1}-x_i) \dfrac{d\overline \bu_{i+1/2}}{dt}+ \bbf (\bu_{i+1})-\bbf(\bu_i)=0, \qquad \dfrac{d\bv_i}{dt}+\Phi_i^{[x_i,x_{i+1}]}+\Phi_i^{[x_{i-1},x_{i}]}=0
 $$
 where $\Phi_i^{[x_i,x_{i+1}]}$ (resp. $\Phi_i^{[x_{i-1},x_{i}]}$) is a consistant approximation of $J^-(\bu_{i} )\dpar{\bv}{x}$ (resp. $J^+(\bu_{i} )\dpar{\bv}{x}$). To define these approximations, we use the data $\bv_{l}$ and $\bv_{l+1/2}$ and the approxiamtons contain in \cite{iserle}. In \cite{AF1} we show that under assumptions similar to the standard Lax Wendroff theorem, that this scheme will provide a sequence that will converge to a weak solution. It is unclear if this scheme can be written in flux form, simply because the volumes associated to the average \remi{already} cover the whole computational domain.
 
\begin{figure}[h]
\subfigure[]{\includegraphics[width=0.5\textwidth]{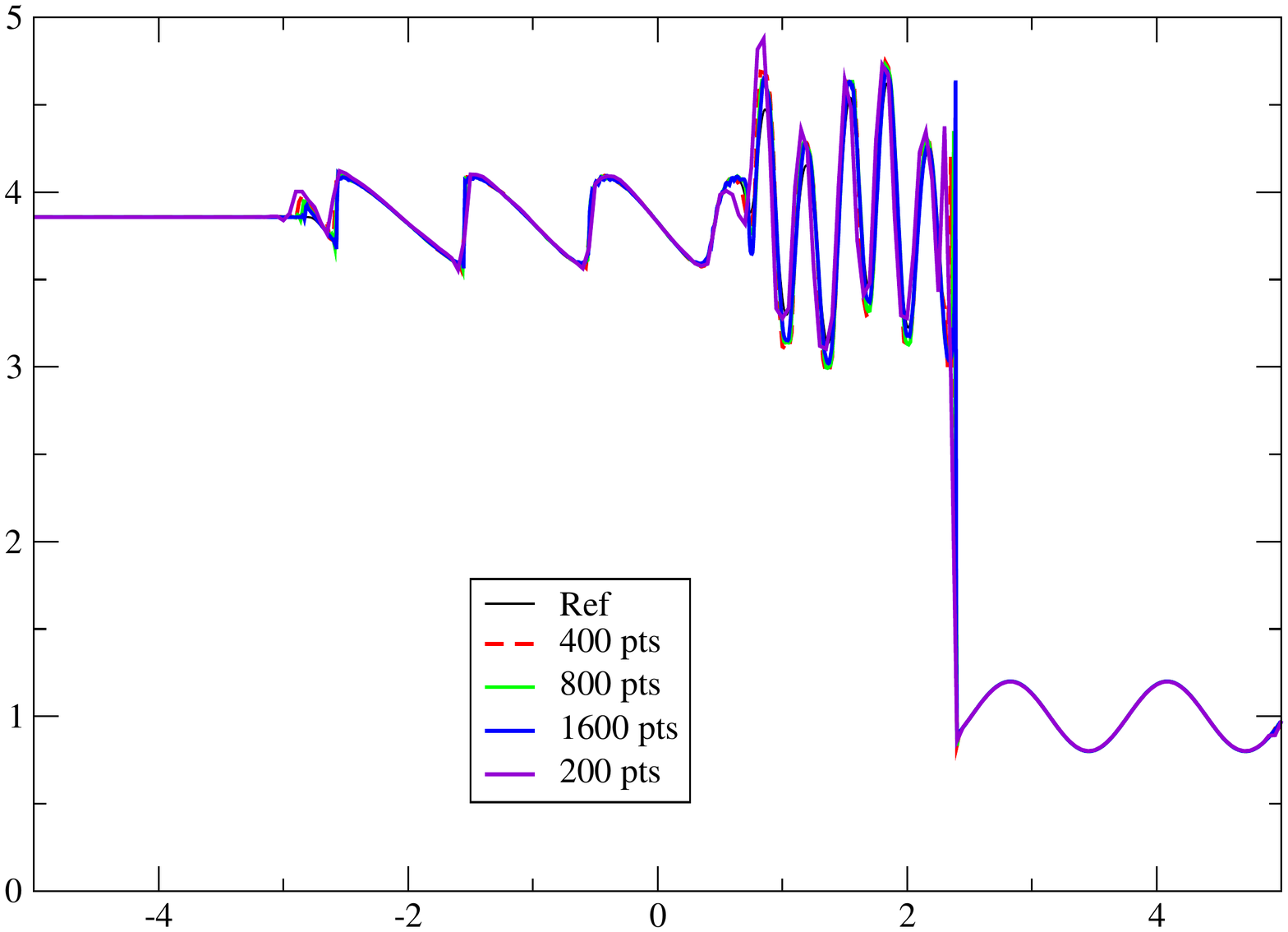}}
\subfigure[]{\includegraphics[width=0.5\textwidth]{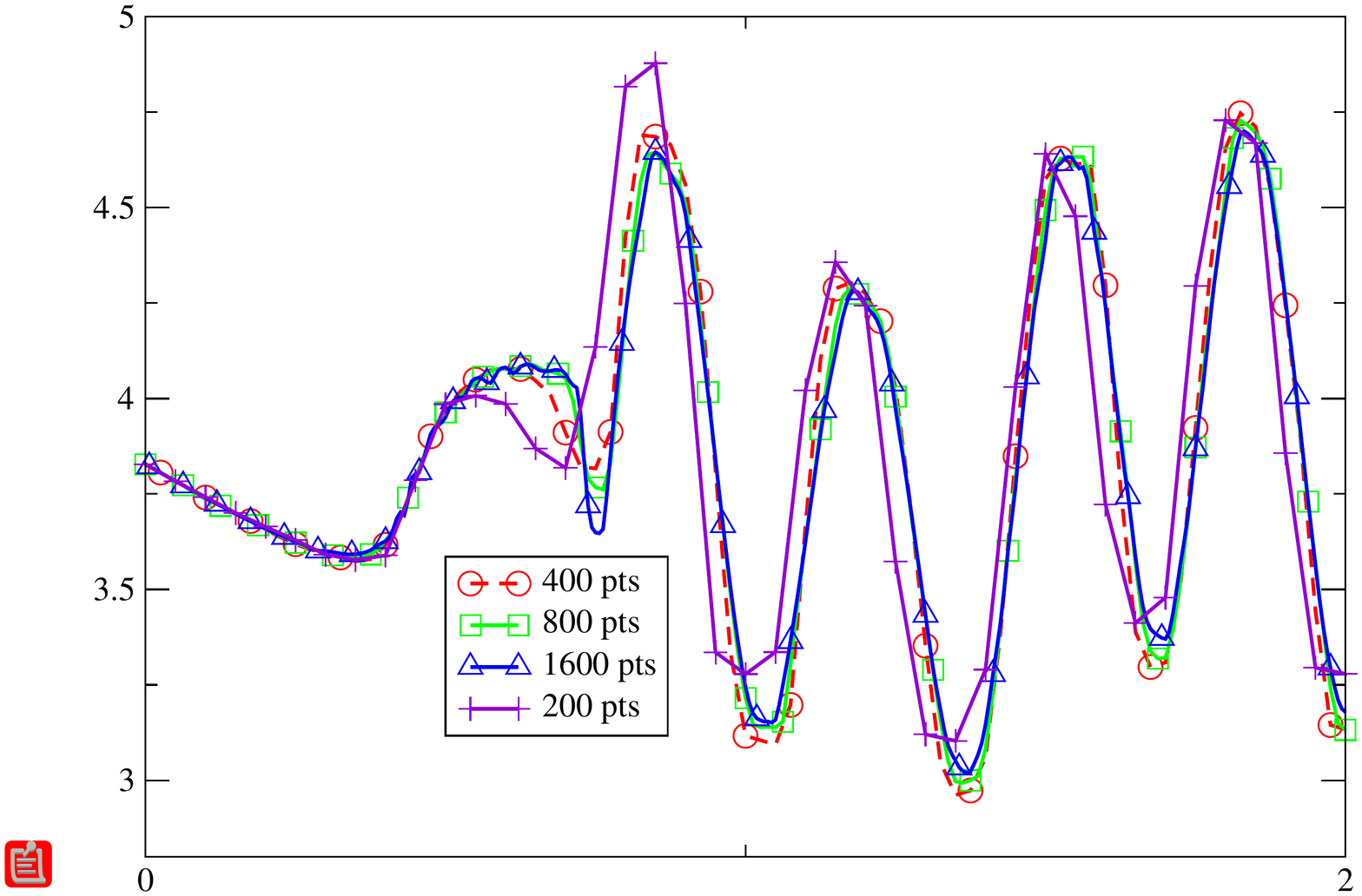}}
\caption{\label{figshu} (a): Solution of the Shu Osher problem , (b): zoom of the solution around the shock.}
\end{figure}

 One illustration, see figure \ref{figshu} is given by the Shu Osher problem where the initial condition are:
$$(\rho, u, p)=\left \{\begin{array}{ll}
(3.857143, 2.629369, 10.3333333) &\text{if } x<-4\\
(1+0.2\sin(5x), 0, 1) &\text{else}
\end{array}\right .
$$
on the domain $[-5,5]$ until $T=1.8$. 
\section*{Acknowledgements} I take this opportunity to thanks my many collaborators and students: M. Ricchiuto (Inria), P.H. Maire (Cea), F. Vilar (Montpellier), R. Loub\`ere (Bordeaux), Ph. \"Offner (Mainz), W. Barsukow (Bordeaux), M. Dumbser (Trento), S. Busto (Vig\`o), S. Tokareva (Los Alamos), P. Baccigaluppi (Milano), L. Micalizzi (Z\"urich), A. Larat (Grenoble), M. Mezine (Numeca), M. Lukakova (Mainz). Without them nothing would have been possible.
\bibliographystyle{plain}
\bibliography{biblio}
\end{document}